\documentclass{amsart}
\usepackage{calc,amssymb,amsthm,amsmath,amsfonts,setspace,mathtools}
\RequirePackage[dvipsnames,usenames]{xcolor}
\usepackage{textgreek}
\usepackage[utf8x]{inputenc}
\DeclareUnicodeCharacter{957}{\textnu}

\usepackage{microtype}

\usepackage{hyperref} \hypersetup{ bookmarks, bookmarksdepth=3,
   bookmarksopen, bookmarksnumbered, pdfstartview=FitH,
   colorlinks,backref,hyperindex, linkcolor=Sepia,
   anchorcolor=BurntOrange, citecolor=MidnightBlue,
   citecolor=OliveGreen, filecolor=BlueViolet, menucolor=Yellow,
   urlcolor=OliveGreen } \usepackage{xspace}
\usepackage{mabliautoref}
\usepackage{colonequals}
\usepackage{verbatim}
\usepackage{enumerate}

\usepackage{fancyvrb}
\DefineVerbatimEnvironment%
{MyVerbatim}{Verbatim}
{formatcom=\color{Violet}}

\usepackage{fancyhdr}

\makeatletter
\def\@settitle{\begin{center}%
  \baselineskip14\p@\relax
    \bfseries
    \normalfont\LARGE
\@title
  \end{center}%
}
\makeatother

\newcommand{\ft}{\operatorname{c}}
\newcommand{\idealm}{\mathfrak{m}}
\DeclareMathOperator{\fpt}{fpt}

\usepackage[textwidth=3.3 cm,textsize=small,shadow
]{todonotes}

\begin{document}
\title[]{The \emph{FrobeniusThresholds} package for \emph{Macaulay2}}

\author[]{Daniel J.\ Hern\'andez}
\address{Department of Mathematics, University of Kansas, Lawrence, KS~66045, USA}
\email{hernandez@ku.edu}
\thanks{D.~J.~Hern\'andez was partially supported by NSF DMS \#1600702.}



\author[]{Karl Schwede}
\address{Department of Mathematics, University of Utah, Salt Lake City, UT~84112, USA}
\thanks{K.~Schwede was supported by NSF CAREER Grant DMS \#1252860/1501102, NSF FRG Grant DMS \#1265261/1501115, NSF grant \#1801849 and a Sloan Fellowship.}
\email{schwede@math.utah.edu}

\author[]{Pedro Teixeira}
\address{Department of Mathematics, Knox College, Galesburg, IL~61401, USA}
\email{pteixeir@knox.edu}

\author[]{Emily E.\ Witt}
\address{Department of Mathematics, University of Kansas, Lawrence, KS~66045, USA}
\email{witt@ku.edu}
\thanks{E.~E.~Witt was partially supported by NSF DMS \#1623035.}
\date{\today}

\begin{abstract}
   This article describes the \emph{Macaulay2} package \emph{FrobeniusThresholds}, designed to estimate and calculate $F$-pure thresholds, more general $F$-thresholds, and related numerical invariants arising in the study of singularities in prime characteristic commutative algebra.
\end{abstract}



\maketitle

\section{Introduction}

This paper describes the \emph{Macaulay2} package \emph{FrobeniusThresholds} \cite{M2, FThresholdsPackage}, which provides tools for computing or estimating certain fundamental invariants in positive characteristic commutative algebra, namely \emph{$F$-pure thresholds}, \emph{$F$-thresholds}, and \emph{$F$-jumping exponents}.
Recall that a ring of prime characteristic $p>0$ is \emph{$F$-pure} if the Frobenius map|that is, the ring endomorphism sending an element to its $p$-th power|is a pure morphism;  under natural geometric hypotheses, this is equivalent to the condition that the Frobenius morphism splits as a map of rings.
The concept of $F$-purity has proven to be important in commutative algebra, and has a rich history.
Indeed, it first appeared in the work of Hochster and Roberts to study local cohomology \cite{HochsterRobertsFrobeniusLocalCohomology}, was compared with rational singularities in  \cite{FedderFPureRat}, and was used to study global properties of Schubert varieties in \cite{MehtaRamanathanFrobeniusSplittingAndCohomologyVanishing}.

After the advent of tight closure \cite{HochsterHunekeTC1}, the use of the Frobenius map to quantify singularity|that is, deviation from regularity|proliferated, and based on a connection discovered between $F$-pure and \emph{log canonical} singularities \cite{HaraWatanabeFRegFPure}, the concept of $F$-purity was generalized to the context of \emph{pairs}.
Along these lines, the \emph{$F$-pure threshold} was defined in analogy with the \emph{log canonical threshold} \cite{TakagiWatanabeFPureThresh}, and \emph{$F$-thresholds} were introduced as a natural extension \cite{MustataTakagiWatanabeFThresholdsAndBernsteinSato}.

The connection between the $F$-pure threshold and the log canonical threshold, however, extends beyond mere analogy.
For example, suppose $h$ is a polynomial with integer coefficients, and that $h_p$ is the polynomial obtained by reducing the coefficients of $h$ modulo a prime $p$.  Then, the $F$-pure thresholds of the reductions $h_p$ converge to the log canonical threshold of $h$ as $p$ tends to infinity \cite{HaraYoshidaGeneralizationOfTightClosure}.
A related result is that of Zhu, who proved that the log canonical threshold of $h$ is at least the $F$-pure threshold of any reduction $h_p$ \cite[Corollary 4.2]{ZhuLogCanoincalThresholdsInPositiveChar}.

This latter result is interesting from a computational perspective, in that it provides lower bounds for log canonical thresholds.
Though there is a general purpose implementation of an algorithm for computing log canonical thresholds in the \emph{Dmodules} package \cite{DmodulesSource}\footnote{The \emph{MultiplierIdeals} package \cite{MultiplierIdealsPackage, MultiplierIdealsArticle} also computes log canonical thresholds in many special cases, including monomial ideals, hyperplane arrangements, generic determinantal ideals, and certain binomial ideals.
}, the function for computing $F$-pure thresholds contained in the \emph{FrobeniusThresholds} package is typically much faster, especially so in low characteristic.

In summary,  the $F$-pure threshold is an interesting numerical invariant, related to many other measures of singularity across all characteristics, and has been the focus of intense study over the past fifteen years.
Unfortunately, it is typically difficult to calculate.
The package \emph{FrobeniusThresholds} is centered on calculating and estimating the $F$-pure threshold and other $F$-thresholds, with the function \texttt{fpt} at its core.
It builds heavily upon the \emph{TestIdeals} package for \emph{Macaulay2} \cite{TestIdealsPackage, TestIdealsPaper}, which provides a broad range of functionality for effective computation in prime characteristic commutative algebra.

\subsection{Some background and notation}
Although some of the functionality implemented in \emph{FrobeniusThresholds} is not restricted to regular ambient rings (see Section~\ref{sec.IsFPT}), for the sake of concreteness, in this introduction  we will work in a polynomial ring over a finite field of  characteristic $p$.
The ideal of this ring generated by its variables is denoted $\idealm$.

Let us outline a way in which natural numerical invariants in prime characteristic commutative algebra are often constructed:  For every natural number $e$, associate to some fixed data---often, a collection of polynomials or ideals---an integer describing something of relevance that depends on $e$ (e.g., the dimension of some interesting vector space constructed in terms of the initial data).
Normalize this integer by dividing by some power of $p^e$, and then take the limit as the integer $e$ tends to infinity.
The resulting limit, if it exists, should then encode some interesting information about the initial data.

Consider a nonzero polynomial $f$ and a natural number $e$.  If $f$ does not vanish at the origin, then set $\nu_f^{\idealm}(p^e) := \infty$.  Otherwise, $f \in \idealm$, and we instead define
\[ \nu_f^{\idealm}(p^e) \coloneqq  \max \{ n \in \mathbb{N} : f^n \notin \idealm^{[p^e]} \}, \]
where $\idealm^{[p^e]}$ denotes the $p^e$-th Frobenius power of $\idealm$, that is, the ideal generated by the $p^e$-th powers of the elements of $\idealm$.
Following our outline, we define
\[ \ft^{\idealm}(f) \coloneqq  \lim_{e \to \infty} \frac{ \nu_f^{\idealm}(p^e)}{p^e}. \]

This limit exists, and is a rational number when $f\in \idealm$, though the latter is far from obvious  \cite{BlickleMustataSmithDiscretenessAndRationalityOfFThresholds}.
Inspired by its connections with the $F$-purity of pairs,  this limit is called the \emph{$F$-pure threshold} of $f$ at the origin.

The $F$-pure threshold is closely related to many other fundamental concepts in prime characteristic commutative algebra.   For instance,
\[ \ft^{\idealm}(f) = \inf \{ t>0 : \tau(f^t)  \subseteq \idealm \} = \sup \{ t>0 : \sigma( f^t) \neq 0 \}, \]
where $\tau(f^t)$ and $\sigma(f^t)$ are the \emph{test ideal} and \emph{$F$-signature}, respectively, associated to $f$ and the formal nonnegative real exponent $t$.  The former is an ideal in the ambient ring of $f$, and the latter is a real number;  both depend on the parameter $t$ and the characteristic $p$ in subtle ways \cite{BlickleMustataSmithDiscretenessAndRationalityOfFThresholds, BlickleSchwedeTuckerFSigPairs1}.

In the literature, the $F$-pure threshold $\ft^{\idealm}(f)$ is often denoted $\fpt(f)$, for obvious reasons.
However, in this note we adopt the former notation to avoid any possible confusion with the function \texttt{fpt} described in Section~\ref{sec.FPT}, which some times does not output the number $\ft^{\idealm}(f) = \fpt(f)$, but returns, instead, lower and upper bounds for that number.

It turns out that the sequence $\big(\nu_f^{\idealm}(p^e) \big)_{e=0}^{\infty}$ itself, and not just its limit, encodes interesting information about $f$.  For example, it is closely related to the \emph{Bernstein--Sato polynomial} of $f$ \cite{MustataTakagiWatanabeFThresholdsAndBernsteinSato}.  Thus, it is natural to ask whether one can recover the sequence $\big( \nu_f^{\idealm}(p^e) \big)_{e=0}^{\infty}$ from the limit $\ft^{\idealm}(f)$.  Fortunately, the answer to this question is yes \cite{MustataTakagiWatanabeFThresholdsAndBernsteinSato, HernandezFPurityOfHypersurfaces}; for each $e$ we have
\begin{equation*}
\nu_f^{\idealm}(p^e) = \lceil p^e \cdot \ft^{\idealm}(f) \rceil - 1.
\end{equation*}

We conclude this subsection by briefly reviewing some natural generalizations.
Suppose that $I$ and $J$ are ideals.
If $I$ is contained in the radical of $J$, then we set
\[ \nu_I^J(p^e) \coloneqq \max \{ n \in \mathbb{N}: I^n \not\subseteq J^{[p^e]} \}, \]
or $\nu_I^J(p^e) \coloneqq 0$, when the set on the right-hand side is empty.
Otherwise, we set $\nu_I^J(p^e) \coloneqq \infty$.
This clearly generalizes the quantity $\nu_f^{\idealm}(p^e)$ considered earlier, and we call
\[ \ft^J(I) \coloneqq  \lim_{e \to \infty} \frac{ \nu_I^J(p^e)}{p^e} \]
the \emph{$F$-threshold of $I$ with respect to $J$}. 
The value $\ft^{\idealm}(I)$ is called the \emph{$F$-pure threshold of $I$} at the origin, and if $I = \langle f \rangle$ is principal,  $\ft^{J}(f) \coloneqq  \ft^J(I)$ is called the \emph{$F$-threshold of $f$ with respect to $J$}.
Like $F$-pure thresholds, $F$-thresholds are rational (when finite), and can be characterized in terms of test ideals.

\subsection*{Acknowledgements.}  The authors enthusiastically thank everyone who helped complete the \emph{FrobeniusThresholds} package: the package coauthors Juliette Bruce and Daniel Smolkin, and contributors  Erin Bela, Zhibek Kadyrsizova, Moty Katzman, Sara Malec, and Marcus Robinson.

Thanks go to the organizers of the \emph{Macaulay2} workshops where much of the functionality described herein was developed, hosted by Wake Forest University in 2012, the University of California, Berkeley in 2014 and 2017, Boise State University in 2015, and the University of Utah in 2016.

Finally, the authors are grateful to the University of Utah for hosting a collaborative development visit in 2018, and to the Institute of Mathematics and its Applications for its generous support for our 2019 Coding Sprint.
The current version of the package was finalized during these events.

\section{The {\tt nu} function}
\label{sec.Nu}

We first consider the \texttt{nu} function, a fundamental component of the package \emph{FrobeniusThresholds}.
We adopt the setup established in the introduction: we work in a polynomial ring $R$ over a finite field of characteristic $p>0$,  $\idealm$ denotes the ideal generated by the variables, and $e$ is a natural number.
If $I$ and $J$ are ideals of $R$, the command \texttt{nu(e,I,J)} outputs the extended integer $\nu_I^J(p^e)$ defined in the introduction; if $f$ is an element of $R$, \texttt{nu(e,f,J)} outputs $\nu_f^J(p^e) \coloneqq \nu_{\langle f \rangle}^J(p^e)$.
When the third argument is omitted from the function \texttt{nu}, it is assumed to be the maximal ideal $\idealm$.

\smallskip
{\small
\setstretch{.67}
\begin{MyVerbatim}
i1 : R = ZZ/11[x,y];

i2 : I = ideal(x^2 + y^3, x*y);

o2 : Ideal of R

i3 : J = ideal(x^2, y^3);

o3 : Ideal of R

i4 : nu(2, I, J)

o4 = 281

i5 : f = x*y*(x^2 + y^2);

i6 : nu(2, f, J)

o6 = 120
\end{MyVerbatim}
}

\subsection{Options for \texttt{nu}}

We now describe the options available for the function \texttt{nu}.
As pointed out in the introduction, if $f \in \mathfrak{m}$, the values $\nu^{\mathfrak{m}}_f(p^e)$ can be recovered from the $F$-pure threshold of $f$.
This is used to speed up computations for certain polynomials whose $F$-pure thresholds can be computed quickly via specialized algorithms or formulas, namely diagonal polynomials, binomials, forms in two variables, and products of factors in simple normal crossing (see Section~\ref{sec.FPT}).
This feature can be disabled by setting the option \texttt{UseSpecialAlgorithms} (default value \texttt{true}) to \texttt{false}.\footnote{In Section~\ref{ss: UseSpecialAlgorithms} we discuss a couple of situations in which this may be desirable.}

\smallskip
{\small
\setstretch{.67}
\begin{MyVerbatim}
i7 : R = ZZ/17[x,y,z];

i8 : f = x^3 + y^4 + z^5; -- a diagonal polynomial

i9 : time nu(10, f)
     -- used 0.0161622 seconds

o9 = 1541642394460

i10 : time nu(10, f, UseSpecialAlgorithms => false)
     -- used 2.06877 seconds

o10 = 1541642394460
\end{MyVerbatim}
}
\smallskip

In general, the function \texttt{nu} works by searching through a list of integers $n$, and checking containments of the $n$-th power of $I$ in the specified Frobenius power of $J$.
It is well known that, for any positive integer $e$,
\[ \nu_I^J(p^e) = \nu_I^J(p^{e-1})\cdot p + L,\]
where the error term $L$ is nonnegative and can be explicitly bounded from above in terms of $p$ and the number of generators of $I$ and $J$.
For instance, the error term $L$ is at most $p-1$ when $I$ is principal and $J$ is arbitrary.
This implies that when searching for the maximal exponent defining \texttt{nu(e,I,J)} for positive $e$, it is safe to start at $p$ times the output of \texttt{nu(e-1,I,J)}, and one need not search too far past this number.
This also suggests that the most efficient way to compute \texttt{nu(e,I,J)} is to compute, successively, \texttt{nu(s,I,J)}, for each integer $s = 0,\ldots,e$, and this is indeed how the computation is done in most cases.

The user can specify how the search is approached through the option \texttt{Search}, which can take two values: 
 \texttt{Binary} (the default value) and \texttt{Linear}.

\smallskip
{\small
\setstretch{.67}
\begin{MyVerbatim}
i11 : R = ZZ/5[x,y,z];

i12 : m = ideal(x, y, z);

o12 : Ideal of R

i13 : time nu(2, m, m^2) -- uses binary search (default)
     -- used 1.82479 seconds

o13 = 97

i14 : time nu(2, m, m^2, Search => Linear) -- but linear search gets luckier
     -- used 0.597035 seconds

o14 = 97
\end{MyVerbatim}
}
\smallskip

If the option \texttt{ReturnList} is changed from its default value of \texttt{false} to \texttt{true}, \texttt{nu} outputs a list of the values $\nu_I^J(p^s)$, for $s=0,\ldots,e$, at no additional computational cost.

\smallskip
{\small
\setstretch{.67}
\begin{MyVerbatim}
i15 : nu(5, x^2*y^4 + y^2*z^7 + z^2*x^8, ReturnList => true)

o15 = {0, 1, 8, 44, 224, 1124}

o15 : List
\end{MyVerbatim}
}
\smallskip

\noindent The same information can be found by setting the option \texttt{Verbose} to \texttt{true}, to request that the values $\nu_I^J(p^s)$ be printed as they are iteratively computed (serving also as a way to monitor the progress of the computation).

As described in the introduction, the integer $\nu_I^J(p^e)$ is the maximal integer $n$ such that the $n$-th power of $I$ is not contained in the $p^e$-th Frobenius power of $J$.  However,
\begin{equation*}
  I^n \subseteq J^{[p^e]} \Longleftrightarrow (I^n)^{[1/p^e]} \subseteq J,
\end{equation*}
where $(I^n)^{[1/p^e]}$ denotes the $p^e$-th \emph{Frobenius root} of $I^n$, as defined in \cite{BlickleMustataSmithDiscretenessAndRationalityOfFThresholds}.
The option \texttt{ContainmentTest} for \texttt{nu} allows the user to choose which of the two types of containment statements appearing above to use toward the calculation of $\nu_I^J(p^e)$.

If \texttt{ContainmentTest} is set to \texttt{StandardPower}, then \texttt{nu(e,I,J)} is computed by testing the left-hand containment above, and when it is set to \texttt{FrobeniusRoot}, the right-hand containment is checked.
For efficiency reasons, the default value for \texttt{ContainmentTest} is set to  \texttt{FrobeniusRoot} if the second argument passed to \texttt{nu} is a polynomial, and is set to \texttt{StandardPower} if the second argument is an ideal.

\smallskip
{\small
\setstretch{.67}
\begin{MyVerbatim}
i16 : R = ZZ/11[x,y,z];

i17 : f = x^3 + y^3 + z^3 + x*y*z;

i18 : time nu(3, f) -- ContainmentTest is set to FrobeniusRoot, by default
     -- used 0.153691 seconds

o18 = 1209

i19 : time nu(3, f, ContainmentTest => StandardPower)
     -- used 10.1343 seconds

o19 = 1209
\end{MyVerbatim}
}

The option \texttt{ContainmentTest} has a third possible value, \texttt{FrobeniusPower}, which allows \texttt{nu} to compute a different but analogous invariant.
The first, third, and fourth authors introduced the notion of a (generalized) Frobenius power $I^{[n]}$ of an ideal $I$, when $n$ is an arbitrary nonnegative integer \cite{hernandez+etal.frobenius_powers}.
When \texttt{ContainmentTest} is set to \texttt{FrobeniusPower}, rather than  computing $\nu_I^J(p^e)$, the function \texttt{nu} computes the maximal integer $n$ for which $I^{[n]}$ is not contained in $J^{[p^e]}$.  This number is denoted $\mu_I^J(p^e)$ in \emph{loc.~cit.}, and equals $\nu_I^J(p^e)$ when $I$ is a principal ideal.
However, these numbers need not agree in general, as we see below.

\smallskip
{\small
\setstretch{.67}
\begin{MyVerbatim}
i20 : R = ZZ/3[x,y];

i21 : m = ideal(x, y);

o21 : Ideal of R

i22 : nu(4, m^5)

o22 = 32

i23 : nu(4, m^5, ContainmentTest => FrobeniusPower)

o23 = 26
\end{MyVerbatim}
}
\smallskip

The last option we describe for \texttt{nu} is \texttt{AtOrigin}.
Recall that $\nu_I^\idealm(p^e)$ can be interpreted as the maximum integer $n$ for which $(I^n)^{[1/p^e]}$ is not contained in $\idealm$. 
When the option \texttt{AtOrigin} is set to \texttt{false} (from its default value \texttt{true}), the function \texttt{nu} determines, instead, the maximum integer $n$ for which $(I^n)^{[1/p^e]}$ is the unit ideal, which can also be characterized as the minimal integer $\nu_I^\mathfrak{n}(p^e)$ as $\mathfrak{n}$ varies among all maximal ideals of the ring.  

\smallskip
{\small
\setstretch{.67}
\begin{MyVerbatim}
i24 : R = ZZ/7[x,y];

i25 : f = (x - 1)^3 - (y - 2)^2;

i26 : nu(3, f)

o26 = infinity

o26 : InfiniteNumber

i27 : nu(3, f, AtOrigin => false)

o27 = 285
\end{MyVerbatim}
}

\section{{\tt isFPT}, {\tt compareFPT} and {\tt isFJumpingExponent}}
\label{sec.IsFPT}

The \emph{FrobeniusThresholds} package contains methods to test candidate values for $F$-pure thresholds and $F$-jumping numbers, even in some singular rings.
Consider a  $\mathbb{Q}$-Gorenstein ring $R$ of characteristic $p>0$, whose index is not divisible by $p$.
Given a parameter $t\in \mathbb{Q}$ and an element $f$ of $R$, the command \texttt{isFPT(t, f)} checks whether $t$ is the $F$-pure threshold of $f$, while \texttt{compareFPT(t, f)} provides further information, returning {\tt-1}, \texttt{0}, or \texttt{1} when $t$ is, respectively, less than, equal to, or greater than the $F$-pure threshold of $f$.  Setting the option \texttt{AtOrigin} to \texttt{true} tells the function to consider the $F$-pure threshold at the origin.  

\smallskip
{\small
\setstretch{.67}
\begin{MyVerbatim}
i1 : R = ZZ/11[x,y,z]/(x^2 - y*(z - 1));

i2 : compareFPT(5/11, z - 1)

o2 = -1

i3 : isFPT(1/2, z - 1)

o3 = true

i4 : isFPT(1/2, z - 1, AtOrigin => true)

o4 = false
\end{MyVerbatim}
}
\smallskip

The general method applied calls upon functionality from the \emph{TestIdeals} package.
The functions \texttt{testIdeal} and \texttt{FPureModule} therein are used to compute the test ideals $\tau(f^t)$ and $\tau(f^{t-\varepsilon})$, for $0<\varepsilon\ll 1$.
By comparing these test ideals, we can decide whether $t$ is the $F$-pure threshold of $f$.
For instance, when $R$ is a polynomial ring, $t$ is the $F$-pure threshold of $f$ at the origin if and only if $\tau(f^t)$ is contained in the homogeneous maximal ideal, but $\tau(f^{t-\varepsilon})$ is not.

Since not only the $F$-pure thresholds, but also the higher $F$-jumping numbers, are determined by containment conditions on test ideals, the functionality is extended to determine whether a given number is an $F$-jumping number.

\smallskip
{\small
\setstretch{.67}
\begin{MyVerbatim}
i5 : R = ZZ/13[x,y];

i6 : f = y*((y + 1) - (x - 1)^2)*(x - 2)*(x + y - 2);

i7 : isFJumpingExponent(3/4, f)

o7 = true

i8 : isFPT(3/4, f)

o8 = false
\end{MyVerbatim}
}

\section{The {\tt fpt} function}
\label{sec.FPT}

The core function in the package \emph{FrobeniusThresholds} is called \texttt{fpt}.  Throughout this section, let $f$ be a polynomial with coefficients in a finite field of characteristic $p$. When passed the polynomial $f$, the function \texttt{fpt} attempts to find the exact value for the $F$-pure threshold of $f$ at the origin, and returns that value, if possible.  Otherwise, it returns lower and upper bounds for the $F$-pure threshold, as demonstrated below.

\smallskip
{\small
\setstretch{.67}
\begin{MyVerbatim}
i1 : R = ZZ/5[x,y,z];

i2 : fpt(x^3 + y^3 + z^3 + x*y*z)

     4
o2 = -
     5

o2 : QQ

i3 : fpt(x^5 + y^6 + z^7 + (x*y*z)^3)

       7  2
o3 = {--, -}
      25  5

o3 : List
\end{MyVerbatim}
}

\subsection{The option \texttt{UseSpecialAlgorithms}}\label{ss: UseSpecialAlgorithms}

The \texttt{fpt} function has an option called \texttt{UseSpecialAlgorithms}, which, when set to \texttt{true} (its default value), tells \texttt{fpt} to first check whether $f$ is a diagonal polynomial, a binomial, a form in two variables, or a product of factors in simple normal crossing, in that order.
When $f$ is a diagonal polynomial, a binomial, or a form in two variables, algorithms of Hern\'andez \cite{HernandezFInvariantsOfDiagonalHyp, HernandezFPureThresholdOfBinomial} or Hern\'andez and Teixeira \cite{HernandezTeixeiraFThresholdFunctions} are executed to compute the $F$-pure threshold.

\smallskip
{\small
\setstretch{.67}
\begin{MyVerbatim}
i4 : fpt(x^17 + y^20 + z^24) -- a diagonal polynomial

      94
o4 = ---
     625

o4 : QQ

i5 : fpt(x^2*y^6*z^10 + x^10*y^5*z^3) -- a binomial

      997
o5 = ----
     6250

o5 : QQ

i6 : R = ZZ/5[x,y];

i7 : fpt(x^2*y^6*(x + y)^9*(x + 3*y)^10) -- a form in two variables

      5787
o7 = -----
     78125

o7 : QQ
\end{MyVerbatim}
}
\smallskip

\noindent The algorithm for computing the $F$-pure threshold of a binary form $f$ requires factoring $f$ into linear forms, and that may be difficult or impossible when that factorization occurs in a Galois field of excessively large order.
This is a situation when the user will want to set the option \texttt{UseSpecialAlgorithms} to \texttt{false}.
However, when a factorization is already known, instead of passing $f$ to \texttt{fpt}, the user can pass a list containing all the pairwise coprime linear factors of $f$ to \texttt{fpt}, and a list containing their respective multiplicities.

\smallskip
{\small
\setstretch{.67}
\begin{MyVerbatim}
i8 : fpt({x, y, x + y, x + 3*y}, {2, 6, 9, 10}) == oo

o8 = true
\end{MyVerbatim}
}
\smallskip

If \texttt{UseSpecialAlgorithms} is set to \texttt{true} and $f$ does not fall into any of the aforementioned cases, then the function \texttt{fpt} next calls \texttt{isSimpleNormalCrossing(f)} (see Section~\ref{subsec.SNC}) to check whether the polynomial $f$ is (locally, at the origin) a product of factors that are in simple normal crossing, in which case the $F$-pure threshold is simply the reciprocal of the largest multiplicity occurring in that factorization.  
Note that the function \texttt{factor} is called whenever \texttt{isSimpleNormalCrossing} is used, and that can sometimes make the verification slow.  The user can avoid this by setting \texttt{UseSpecialAlgorithms} to \texttt{false}.

\subsection{When no special algorithm applies}

We now explain how the function  \texttt{fpt} proceeds when no special algorithm is available, or when \texttt{UseSpecialAlgorithms} is set to \texttt{false}.
In this case, \texttt{fpt} computes a sequence of lower and upper bounds for the $F$-pure threshold of $f$, and either finds its exact value in this process, or outputs the last of these sets of bounds, which will be the tightest among all computed.
The value of the option \texttt{DepthOfSearch} determines the precision of the initial set of bounds, and the option \texttt{Attempts} determines, roughly, how many new, tighter sets of bounds are to be computed.

More specifically, let $e$ denote the value of the option \texttt{DepthOfSearch}, which conservatively defaults to \texttt{1}.
The \texttt{fpt} function first computes $\nu=\nu_f(p^e)$, which agrees with the output of \texttt{nu(e,f)}.
It is well known that the $F$-pure threshold of $f$ is greater than $\nu/p^e$ and at most $(\nu+1)/p^e$, and applying  \cite[Proposition~4.2]{HernandezFPurityOfHypersurfaces} to the lower bound tells us that the $F$-pure threshold of $f$ must be at least $\nu/(p^e-1)$.
In summary, we know that the $F$-pure threshold of $f$ must lie in the closed interval
\begin{equation}
\label{estimating-interval: e}
\tag{$\dagger$}
\left[ \frac{\nu}{p^e-1}, \frac{\nu+1}{p^e} \right].
\end{equation}

With these estimates in hand, the subroutine \texttt{guessFPT} is called to make some ``educated guesses" in an attempt to identify the $F$-pure threshold within this interval, or at least narrow down this interval to produce improved estimates.  The number of ``guesses" is controlled by the option \texttt{Attempts}, which conservatively defaults to \texttt{3}.  If \texttt{Attempts} is set to \texttt{0}, then \texttt{guessFPT} is bypassed. If  \texttt{Attempts} is set to at least \texttt{1}, then a first check is run to verify whether the right-hand endpoint $(\nu+1)/p^e$ of the above interval \eqref{estimating-interval: e} is the $F$-pure threshold.  We illustrate this below.

\smallskip
{\small
\setstretch{.67}
\begin{MyVerbatim}
i9 : f = x^2*(x + y)^3*(x + 3*y^2)^5;

i10 : fpt(f, Attempts => 0) -- a bad estimate

          1
o10 = {0, -}
          5

o10 : List

i11 : fpt(f, Attempts => 0, DepthOfSearch => 3) -- a better estimate

        21   22
o11 = {---, ---}
       124  125

o11 : List

i12 : fpt(f, Attempts => 1, DepthOfSearch => 3) -- the right-hand endpoint
      (ν+1)/p^e is the F-pure threshold

       22
o12 = ---
      125

o12 : QQ
\end{MyVerbatim}
}
\smallskip

If  \texttt{Attempts} is set to at least \texttt{2} and the right-hand endpoint $(\nu+1)/p^e$ of the interval \eqref{estimating-interval: e} is not the $F$-pure threshold, then a second check is run to verify whether the left-hand endpoint $\nu/(p^e-1)$ of this interval is the $F$-pure threshold.

\smallskip
{\small
\setstretch{.67}
\begin{MyVerbatim}
i13 : f = x^6*y^4 + x^4*y^9 + (x^2 + y^3)^3;

i14 : fpt(f, Attempts => 1, DepthOfSearch => 3)

       17   7
o14 = {--, --}
       62  25

o14 : List

i15 : fpt(f, Attempts => 2, DepthOfSearch => 3) -- the left-hand endpoint
      ν/(p^e-1) is the F-pure threshold

      17
o15 = --
      62

o15 : QQ
\end{MyVerbatim}
}
\smallskip

If neither endpoint is the $F$-pure threshold and \texttt{Attempts} is set to more than \texttt{2}, then  additional checks are performed at certain numbers within the interval.
First, a number in the interval is selected, according to criteria specified by the value of the option \texttt{GuessStrategy}; we refer the reader to the documentation of this option for more details.
Then the function \texttt{compareFPT} is used to test that number. If that ``guess'' is correct, its value is returned; otherwise, the information returned by \texttt{compareFPT} is used to narrow down the interval, and this process is repeated as many times as specified by \texttt{Attempts}.

\smallskip
{\small
\setstretch{.67}
\begin{MyVerbatim}
i16 : f = x^3*y^11*(x + y)^8*(x^2 + y^3)^8;

i17 : fpt(f, DepthOfSearch => 3, Attempts => 4)

        1   4
o17 = {--, --}
       20  75

o17 : List

i18 : fpt(f, DepthOfSearch => 3, Attempts => 6)

        13   4
o18 = {---, --}
       250  75

o18 : List

i19 : fpt(f, DepthOfSearch => 3, Attempts => 8) 

       1
o19 = --
      19

o19 : QQ
\end{MyVerbatim}
}
\smallskip

The option \texttt{Bounds} allows the user to specify known lower and upper bounds for the $F$-pure threshold of $f$, in order to speed up computations or to refine previously obtained estimates.

\smallskip
{\small
\setstretch{.67}
\begin{MyVerbatim}
i20 : f = x^7*y^5*(x + y)^5*(x^2 + y^3)^4;

i21 : fpt(f, DepthOfSearch => 3, Attempts => 5)

        19   1
o21 = {---, --}
       250  13

o21 : List

i22 : fpt(f, DepthOfSearch => 3, Attempts => 5, Bounds => oo)

        45   1
o22 = {---, --}
       589  13

o22 : List
\end{MyVerbatim}
}
\smallskip

If \texttt{guessFPT} is unsuccessful and \texttt{FinalAttempt} is set to \texttt{true}, the \texttt{fpt} function proceeds to use the convexity of the $F$-signature function and a secant line argument to attempt to narrow down the interval bounding the $F$-pure threshold.
If successful, the new lower bound may coincide with the upper bound, in which case we can conclude that it is the desired $F$-pure threshold.
If this is not the case, a check is performed to verify if the new lower bound is the $F$-pure threshold.

\smallskip
{\small
\setstretch{.67}
\begin{MyVerbatim}
i23 : f = 2*x^10*y^8 + x^4*y^7 - 2*x^3*y^8;

i24 : numeric fpt(f, DepthOfSearch => 3)

o24 = {.14, .144}

o24 : List

i25 : numeric fpt(f, DepthOfSearch => 3, FinalAttempt => true)
      -- FinalAttempt improves the estimate slightly

o25 = {.142067, .144}

o25 : List
\end{MyVerbatim}
}
\smallskip

The computations performed when \texttt{FinalAttempt} is set to \texttt{true} are often slow, and often fail to improve the estimate, and for this reason, this option should be used sparingly.
It is typically more effective to increase the values of \texttt{Attempts} or \texttt{DepthOfSearch}, instead.

\smallskip
{\small
\setstretch{.67}
\begin{MyVerbatim}
i26 : time numeric fpt(f, DepthOfSearch => 3, FinalAttempt => true)
     -- used 0.72874 seconds

o26 = {.142067, .144}

o26 : List

i27 : time fpt(f, DepthOfSearch => 3, Attempts => 7) 
     -- used 0.452872 seconds

      1
o27 = -
      7

o27 : QQ

i28 : time fpt(f, DepthOfSearch => 4) 
     -- used 0.338834 seconds

      1
o28 = -
      7

o28 : QQ
\end{MyVerbatim}
}
\smallskip

As seen in several examples above, when the exact answer is not found, a list containing the endpoints of an interval containing the $F$-pure threshold of $f$ is returned.
Whether that interval is open, closed, or a mixed interval depends on the options passed (it will be \emph{open} whenever \texttt{Attempts} is set to at least \texttt{3}); if the option \texttt{Verbose} is set to true, the precise interval will be printed.

\smallskip
{\small
\setstretch{.67}
\begin{MyVerbatim}
i29 : fpt(f, DepthOfSearch => 3, FinalAttempt => true, Verbose => true)

Starting fpt ...

fpt is not 1 ...

Verifying if special algorithms apply...

Special fpt algorithms were not used ...

ν has been computed: ν = nu(3,f) = 17 ...

fpt lies in the interval [ν/(p^e-1),(ν+1)/p^e] = [17/124,18/125] ...

Starting guessFPT ...

The right-hand endpoint is not the fpt ...

The left-hand endpoint is not the fpt ...

guessFPT narrowed the interval down to (7/50,18/125) ...

Beginning F-signature computation ...

First F-signature computed: s(f,(ν-1)/p^e) = 793/15625 ...

Second F-signature computed: s(f,ν/p^e) = 342/15625 ...

Computed F-signature secant line intercept: 8009/56375 ...

F-signature intercept is an improved lower bound;
Using F-regularity to check if it is the fpt ...

The new lower bound is not the fpt ...

fpt failed to find the exact answer; try increasing the value of
    DepthOfSearch or Attempts.

fpt lies in the interval (8009/56375,18/125).

        8009   18
o29 = {-----, ---}
       56375  125

o29 : List
\end{MyVerbatim}
}
\smallskip

Finally, we point out that one can set the option \texttt{AtOrigin} from its default value of \texttt{true} to \texttt{false}, to compute the $F$-pure threshold globally. In other words, it computes the minimum of the $F$-pure threshold at all maximal ideals.

\smallskip
{\small
\setstretch{.67}
\begin{MyVerbatim}
i30 : R = ZZ/7[x,y];

i31 : f = x*(y - 1)^2 - y*(x - 1)^3;

i32 : fpt(f)

o32 = 1

i33 : fpt(f, AtOrigin => false)

      5
o33 = -
      6

o33 : QQ
\end{MyVerbatim}
}
\smallskip

\noindent In this case, most features enabled by \texttt{UseSpecialAlgorithms => true}
are ignored, except for the check for simple normal crossings;  \texttt{FinalAttempt => true}
is also ignored.

\subsection{The function \texttt{isSimpleNormalCrossings}} \label{subsec.SNC}
As mentioned earlier, the function \texttt{isSimpleNormalCrossings} verifies whether a polynomial $f$ is a product of factors in simple normal crossing.
Suppose that $f$ has factorization $f_1^{a_i} f_2^{a_2} \cdots f_n^{a_n}$.  Recall that its factors $f_i$ are said to be in
\emph{simple normal crossing} at a point if, locally, they form part of a regular system of parameters.  The function \texttt{isSimpleNormalCrossings} determines whether $f$ has simple normal crossings at the origin by computing the Jacobian matrix of each subset of $\{ f_1, \ldots, f_n \}$ (evaluated at the origin), and checking that these matrices have the expected rank, and that these subsets generate ideals of the appropriate height.

\smallskip
{\small
\setstretch{.67}
\begin{MyVerbatim}
i34 : R = ZZ/7[x,y,z];

i35 : isSimpleNormalCrossing(x^2 - y^2)

o35 = true

i36 : isSimpleNormalCrossing(x^2 - y*z)

o36 = false
\end{MyVerbatim}
}
\smallskip

The function \texttt{isSimpleNormalCrossing} is exposed to the user, so can be used independent of any $F$-pure threshold calculation.
If the user sets its option \texttt{AtOrigin} to \texttt{false} (its default value is \texttt{true}), then the function checks whether the $f_i$ are in simple normal crossing \emph{everywhere}, which can be much slower, since Jacobian ideals are computed.

\smallskip
{\small
\setstretch{.67}
\begin{MyVerbatim}
i37 : R = QQ[x,y,z];

i38 : f = (y - (x - 1)^2)*y^2; --SNC at the origin, but not globally

i39 : isSimpleNormalCrossing(f)

o39 = true

i40 : isSimpleNormalCrossing(f, AtOrigin => false)

o40 = false
\end{MyVerbatim}
}
\smallskip

\section{Possible Future Directions}
\label{sec.FutureDirections}

One natural direction of development would be to  incorporate the test ideals $\tau(I^t)$ when computing $F$-thresholds in the case where the ideal $I$ is nonprincipal.
The theoretical foundation for computing such test ideals has already largely been worked out in \cite{SchwedeTuckerTestIdealsOfNonPrincipal}, but such an update to the \emph{FrobeniusThresholds} package would require the  \emph{TestIdeals} package to be updated first.

At the moment, the specialized algorithm called by the \texttt{fpt} function for computing the $F$-pure threshold of a homogeneous polynomial in two variables limits us to working with the standard grading on the ambient polynomial ring.  Clearly, it is desirable to be able to work with nonstandard gradings as well, and results of \cite{HernandezTeixeiraFThresholdFunctions} suggest an algorithm for doing so.  Thus, the package could be improved by finalizing and implementing such an algorithm.

Finally, it would be desirable to develop and implement additional algorithms for computing $F$-pure thresholds and $F$-jumping numbers for additional classes of polynomials.  The first, third, and fourth authors, along with Josep \'Alvarez Montaner, Jack Jeffries, and Luis N\'u\~nez-Betancourt,  are currently working on developing such algorithms.  The theoretical foundation of these algorithms lies in polyhedral geometry and integer programming, making them natural candidates for  implementation in \emph{Macaulay2}.

\bibliographystyle{skalpha}
\bibliography{MainBib}

\end{document}